\documentclass[12pt]{article}
\usepackage{amsmath,amssymb}
\numberwithin{equation}{section}

\begin{document}
\newcommand{\derl}{\partial_L}
\newcommand{\derr}{\partial_R}
\newtheorem{prop}{Proposition}

\def\tr{{\rm tr}\,}
\def\ba{\begin{array}}
\def\ea{\end{array}}
\def\be{\begin{equation}}
\def\ee{\end{equation}}
\def\bea{\begin{eqnarray}}
\def\eea{\end{eqnarray}}
\def\wc{\bar{w}}
\def\ld{\ldots}
\def\ie{i.e.}
\def\C{{\mathbb C}}
\def\Z{{\mathbb Z}}
\def\R{{\mathbb R}}
\def\cf{{\mathcal{F}}}
\def\cs{{\mathcal{S}}}
\def\cl{{\mathcal{L}}}
\def\cx{{\mathcal{X}}}
\def\cz{{\mathcal{Z}}}
\def\su{{\mathfrak s \mathfrak u}}
\def\fd{{f^\dagger}}
\def\dg{{^\dagger}}
\def\fdo{{f^{0\dagger}}}
\def\1{{\bf 1}}

\title{Description of surfaces associated with Grassmannian sigma models on Minkowski space}

\author{
A.~M. Grundland\thanks{email address: grundlan@crm.umontreal.ca}
\\
Centre de Recherches Math{\'e}matiques, Universit{\'e} de Montr{\'e}al, \\
C. P. 6128, Succ.\ Centre-ville, Montr{\'e}al, (QC) H3C 3J7, Canada \\
Universit\'{e} du Qu\'{e}bec, Trois-Rivi\`{e}res CP500 (QC) G9A 5H7, Canada \\
{\rm and} \\
L. \v{S}nobl\thanks{email address: Libor.Snobl@fjfi.cvut.cz}
\\
Centre de Recherches Math{\'e}matiques, Universit{\'e} de Montr{\'e}al, \\
C. P. 6128, Succ.\ Centre-ville, Montr{\'e}al, (QC) H3C 3J7, Canada \\
Faculty of Nuclear Sciences and Physical Engineering, \\
Czech Technical University, \\
B\v rehov\'a 7, 115 19 Prague 1, Czech Republic\\} \date{}

\maketitle

\abstract{We construct and investigate smooth orientable surfaces in $\su(N)$ algebras. 
The structural equations of surfaces associated with Grassmannian sigma models on Minkowski space 
are studied using moving frames adapted to the surfaces. The first and second fundamental forms 
of these surfaces as well as the relations between them as expressed in the Gauss--Weingarten and 
Gauss--Codazzi--Ricci equations are found. The scalar curvature and the mean curvature vector expressed 
in terms of a solution of Grassmanian sigma model are obtained.} 
\smallskip

\smallskip


\section{Introduction}

Sigma models are of great interest in mathematical physics because a significant number of physical systems 
can be reduced to these, relatively simple,  models, either on Euclidean or Minkowski space.  
One such example is the string theory in which 
sigma models on spacetime and their supersymmetric extensions play a crucial role.
Other relevant applications of recent interest are in the areas of statistical physics (for example reduction of self--dual 
Yang--Mills equations to the Ernst model \cite{Abl,David}), phase transitions \cite{Nel,Cha} 
and the theory of fluid membranes \cite{Ou,Saf}. 

The objective of this paper is to study geometric properties of surfaces in Lie algebras associated with sigma models on Minkowski space. 
Recently, we investigated surfaces in $\su(N)$ associated with $\C P^{N-1}$ sigma models \cite{Grusno} and found a few examples
\cite{Grusnosym}. In this paper we extend this approach to more general models based 
on Grassmannian manifolds, i.e. the homogeneous spaces 
$$ G(m,n) = \frac{SU(N)}{S(U(m) \times U(n))}, \; \; N=m+n.$$
Grassmannian sigma models are a generalization of $\C P^{N-1}$ sigma models. Their important common feature is that the Euler--Lagrange
equations can be written in terms of projectors only \cite{Sas}. They share a lot of properties like infinite number of 
local and/or nonlocal conserved quantities, Hamiltonian structure, complete integrability, infinite--dimensional symmetry algebra, 
existence of multisoliton solutions etc. The $N\times N$ projector matrix $P$ for the complex Grassmannian sigma models has 
in general rank lower than the corresponding one for the $\C P^{N-1}$ sigma model and consequently new phenomena can arise.

The generalization of our previous results \cite{Grusno,Grusnosym} to Grassmannian sigma models seemed to be rather natural -- in fact,
it was in a sense more straightforward than the generalization from $\C P^{1}$ to $\C P^{N-1}$, provided one expressed the 
corresponding formulas in terms of the projector (\ref{defP}). On the other hand, a different perspective  obtained in more general 
case allowed to write some of the results in more compact and presumably more natural way.

The results can be of interest in the area of relativistic classical and quantum field theory \cite{Din,Gro}, string theory in which
sigma models on spacetime and their supersymmetric extensions play a crucial role \cite{Bar}. Other relevant applications
of recent interest are in the areas of nonlinear interactions in particle physics \cite{Bre}. 
The explicit forms of the surfaces 
can serve to illuminate the role of the Kac--Moody algebras in integrable models associated with 
the Grassmannian sigma models \cite{Lus,Koi}.

The paper is organized as follows. In Section \ref{SecGrsm} we recall some basic notions and definitions dealing with the 
complex Grassmannian sigma models and their Euler--Lagrange equations. In Section \ref{Secsurf} we perform the analysis
of two--dimensional surfaces immersed in the $\su(N)$ algebra, associated with these models. The geometric properties 
of surfaces and the construction of moving frames are discussed in detail in Sections \ref{SecGW},\ref{Secmf}. Finally, we summarize
our results.

\section{Grassmannian sigma models and their Euler--Lagrange equations}\label{SecGrsm}

As a starting point let us present some basic formulae and notation for complex Grassmannian sigma models defined on Minkowski space.
We adapt to our signature the notation introduced in \cite{Sas} for Euclidean Grassmannian sigma models.

The Grassmannian manifold is defined as homogeneous space
\be\label{Gmn}
G(m,n) = \frac{SU(N)}{S(U(m) \times U(n))}, \; \; N=m+n.
\ee
We express elements $G(m,n)$ using the equivalence classes of elements $g \in SU(N)$ as
\be\label{eqcl}
[g]= \{ g.\psi \; | \; \psi=\left( \ba{cc} U_m & 0 \\ 0 & U_n \ea \right), \; U_m \in U(m), U_n \in U(n), {\det} \psi=1 \}.
\ee
We decompose $g \in SU(N)$ into submatrices $X,Y$
\be
g=(\phi_1,\ldots,\phi_N)=(X,Y), \; \; X=(\phi_1,\ldots,\phi_m), \; Y=(\phi_{m+1},\ldots,\phi_N)
\ee
and from $g^\dagger g = \1$, i.e. $\phi^\dagger_j \phi_k=\delta_{jk}$ we find
$$ X\dg X = \1_{m\times m}, \; X\dg Y=0, \; Y\dg X=0, \; Y\dg Y=\1_{n\times n}.$$
From these orthogonality relations and (\ref{eqcl}) we realize that on the subset of $G(m,n)$ such that
the lower square $n \times n$ submatrix of $Y$ is nonsingular, 
$X$ itself is sufficient to determine $[g]$ (since $U_n$ can be used to bring the lower square part of $Y$
to $\1_{n\times n}$ and the remaining entries in $Y$ are fully determined by the orthogonality properties).
In the following we shall assume that we are working in such chart. Evidently they cover 
the whole $G(m,n)$ up to lower dimensional submanifolds. We shall denote the equivalence classes
either $[X]$ or $[g]$ depending on circumstances.  Note that there is still some freedom in the choice of $X$, namely
$X$ and $X.h$, $h \in \left( \ba{cc} SU(m) & 0 \\ 0 & \1 \ea \right)$ give rise to the same equivalence class $[X]=[Xh]$.
Therefore, one cannot identify $X=[X]$.

Let $\xi^0$, $\xi^1$ be the standard Minkowski coordinates in $\R^2$, with the metric
$$ ({\rm d}s)^2 =({\rm d}\xi^0)^2-({\rm d}\xi^1)^2. $$ 
In what follows we suppose that
$\xi_L=\xi^0+\xi^1 ,\, \xi_R=\xi^0-\xi^1$ are the light--cone coordinates in $\R^2$, i.e.     
\be\label{Minkmetr}
( {\rm d}s)^2= {\rm d}\xi_L {\rm d}\xi_R. 
\ee
We shall denote by $\derl$ and $\derr$ the derivatives with respect to $\xi_L$ and $\xi_R$, respectively.

Let us assume that $\Omega$ is an open, connected and simply connected subset in
$\R^2$ with Minkowski metric (\ref{Minkmetr}).
We define covariant derivatives $D_\mu$ acting on maps $X: \Omega \rightarrow G(m,n)$ by
\be
D_\mu X = \partial_{\mu} X - X X^{\dagger} \partial_\mu X, \qquad \partial_{\mu} \equiv \partial_{\xi^\mu}, \ \mu=0,1.
\ee

In the study of Grassmannian sigma models we are interested in maps $X:\Omega\rightarrow G(m,n)$
 which are stationary points of the action functional 
\be\label{GRmodel}
{\cs} = \int_\Omega \tr \{ (D_\mu X)^{\dagger} (D^\mu X) \} {\rm d}\xi^0 {\rm d}\xi^1.
\ee
The Lagrangian density can be further developed to get
\be\label{Lagr}
{\cal L} = \tr \{ (D_\mu X)^{\dagger} (D^\mu X) \} = \tr \{\partial^\mu X (\partial_\mu X)^{\dagger} P \}
\ee
where 
\be\label{defP}
P = \1-X X\dg
\ee
is an orthogonal projector, i.e. $ P^2=P, \; P^\dagger=P$ satisfying $PX=0, \, X\dg P=0$.

The action (\ref{GRmodel}) has the local (gauge) $SU(m)$ symmetry 
\be\label{Um}
X(\xi_L,\xi_R) \rightarrow X(\xi_L,\xi_R).h(\xi_L,\xi_R), \; h(\xi_L,\xi_R) \in \left( \ba{cc} SU(m) & 0 \\ 0 & \1 \ea \right)
\ee
proving that the model doesn't depend on the choice of representatives $X$ of elements $[X]$ of $G(m,n)$;
and the $SU(N)$ global symmetry
\be\label{SUN}
 X \rightarrow g X, \; g \in SU(N).
\ee 
It is also invariant under the conformal transformations 
\be\label{Poincare}
 \xi_L \rightarrow \alpha (\xi_L), \; \xi_R \rightarrow \beta ( \xi_R ),  
\ee
where $\alpha,\beta: \R \rightarrow \R$ are arbitrary 1--to--1 maps such that $\derl \alpha(\xi_L)\neq 0, \ \derr \beta(\xi_R) \neq 0$,
as well as under the parity transformation 
\be\label{parity}
 \xi_L \rightarrow \xi_R, \; \xi_R \rightarrow \xi_L .
\ee
Let us note that the invariance properties (\ref{Um})--(\ref{parity}) are naturally reproduced on the level of Euler--Lagrange equations.

By variation of the action (\ref{GRmodel}) respecting the constraint
\be\label{constraint}
X\dg X=\1, \; \; {\rm i.e.} \; \; \delta X\dg X+X \dg \delta X=0, \; \partial_{\mu} X\dg X+X \partial_{\mu} \delta X=0
\ee
and assuming that due to suitable boundary conditions the boundary terms vanish we find
the Euler--Lagrange equations 
\be\label{ELeq}
 P ( \derl \derr X -  2 \partial_{\mu} X X\dg \partial^{\mu} X ) = 0.
\ee
They can be also expressed in the matrix form
\be\label{ELeqmat}
[\derl \derr P,P] =0
\ee
or in the form of a conservation law
\be\label{ELeqcons}
\derl [\derr P, P] + \derr [\derl P,P] =0.
\ee
Methods for finding special solutions of (\ref{ELeq}), e.g. 
soliton solutions, are known \cite{Pie,Jac}.

By explicit calculation one can check that the real--valued functions
\be\label{currents}
J_L =  \tr( \derl X \derl X^\dagger P  ), \; \; J_R = \tr( \derr X \derr X^\dagger P )
\ee
satisfy
\be\label{conslaw}
 \derl J_R = \derr J_L =0 
\ee
for any solution $X$ of the Euler--Lagrange equations (\ref{ELeq}). The functions $J_L,J_R$ are invariant under local 
$SU(m)$ and global $SU(N)$ transformations (\ref{Um}) and (\ref{SUN}).

\section{Surfaces obtained from Grassmannian sigma model}\label{Secsurf}

Let us now discuss the analytical description of a two--dimensional smooth orientable surface $\cf$ immersed in the $\su(N)$ algebra, 
associated with the Grassmannian sigma model (\ref{ELeq}).
We shall construct an exact $\su(N)$--valued 1--form whose ``potential'' 0--form defines the surface $\cf$. 
Next, we shall investigate the geometric characteristics of the surface $\cf$. 

Let us introduce a scalar product
$$ ( A , B ) = -\frac{1}{2} \tr AB $$
on $\su(N)$ and identify the $(N^2-1)$--dimensional Euclidean space with the $\su(N)$ algebra
$$ \R^{N^2-1} \simeq \su(N). $$
We denote 
\be\label{mlmr} M_L = [\derl P,P], \qquad M_R = [\derr P,P]. \ee
It follows from (\ref{ELeqcons}) that if $X$ is a solution of the Euler--Lagrange equations (\ref{ELeq}) then
\be\label{ELeqcons2}
\derl M_R+ \derr M_L =0 .
\ee
We identify tangent vectors to the surface $\cf$ with the matrices $M_L$ and $M_R$, as follows
\be\label{dlxdrx}
 {\cz}_L=M_L, \; \; {\cz}_R=-M_R .
\ee
Equation (\ref{ELeqcons2}) implies there exists a closed $\su(N)$--valued 1--form on $\Omega$
$$ {\cz} = {\cz}_L {\rm d} \xi_L + {\cz}_R {\rm d} \xi_R, \; \;  {\rm d} {\cz} =0 .$$
Because ${\cz}$ is closed and $\Omega$ is connected and simply connected, ${\cz}$ is also exact. In other words,
there exists a  well--defined 
$\su(N)$--valued function $Z$ on $\Omega$ such that ${\cz} = {\rm d} Z$. 
The matrix function $Z$ is unique up to addition of any constant element of $\su(N)$ and 
we identify the components of $Z$ with the coordinates of the sought--after surface $\cf$ in $\R^{N^2-1}$. 
Consequently, we get 
\be\label{cldx}
\derl Z= {\cz}_L,\ \derr Z= {\cz}_R.
\ee
The map $Z$ is called the Weierstrass formula for immersion. 
In practice, the surface $\cf$ is found by integration
\be\label{surfaceN}
\cf: \;  Z(\xi_L,\xi_R) = \int_{\gamma(\xi_L,\xi_R)} {\cz}
\ee
along any curve $\gamma(\xi_L,\xi_R)$ in $\Omega$ connecting the point $(\xi_L,\xi_R)\in \Omega$ with an arbitrary chosen
point $(\xi_L^0,\xi_R^0)\in \Omega$.

By computation of traces of ${\cz}_{B}.{\cz}_{D}, \ B,D = L,R$ we find the components of the induced metric on the surface 
${\cf}$
\bea\label{metric} 
 G   & = & \left( \begin{array}{cc} G_{LL}, & G_{LR} \\ G_{LR}, & G_{RR} \end{array}  \right) = \\ 
\nonumber & = & \left( \begin{array}{cc} J_L  & -  \tr \left( \frac{\derl X \derr X\dg +\derr X \derl X\dg}{2} P \right)   \\ 
-  \tr \left( \frac{\derl X \derr X\dg +\derr X \derl X\dg}{2} P \right)  
 & J_R \end{array}  \right). 
\eea
The first fundamental form of the surface $\cf$ takes surprisingly compact form
\bea
\nonumber I   &= &   J_L ({\rm d}\xi_L)^2 -  2 G_{LR} {\rm d}\xi_L {\rm d}\xi_R +  J_R ({\rm d} \xi_R)^2 \\
& = & (2\delta_{B,D}-1) \tr(\partial_{B} X \partial_{D} X\dg P) {\rm d}\xi_B {\rm d}\xi_D
\label{1stff}
\eea
where summation over repeated indices $B,D=L,R$ applies and $\delta_{B,D}=1$ if $B=D$ and 0 otherwise.

In order to establish conditions on a solution $X$ of the Euler--Lagrange equations (\ref{ELeq}) under which 
the surface exists, we introduce a scalar product on the space of $N\times m$ matrices $X$
$$(b,a) = \tr(a . b\dg), \; a,b \in \C^{N\times m}$$
and employ the Schwarz inequality, i.e. 
\be\label{Schw}
 |\tr( a b\dg A)|^2 \leq \tr(a a\dg A)\tr( b b \dg A)
\ee
valid for any positive hermitean operator $A$, namely for $P: P(a)=P.a$.
We may write
\be\label{jdse}
 J_D =\tr(  \partial_D X \partial_D X\dg P ) \geq 0, \; D=L,R 
\ee
and 
\be\label{dgse} 
\det G = \tr( \derl X \derl X\dg P)\tr( \derr X \derr X\dg P) - \left( \Re \ \tr ( \derl X \derr X\dg  P ) \right)^2 \geq 0
\ee
since
$$ \tr( \derl X \derl X\dg P) \ \tr(\derr X \derr X\dg P)\geq | \tr(\derl X \derr X\dg P )|^2 $$
$$ \geq  \left( 
\Re \ \tr ( \derl X \derr X\dg P ) \right)^2.$$
Therefore the first fundamental form $I$ defined by (\ref{1stff}) is positive for any solution $X$ 
of the Euler--Lagrange equations (\ref{ELeq}).

Analyzing the cases when equalities in Schwarz inequality hold we find that $I$ is 
positive definite in the point $(\xi^0_L,\xi^0_R)$ either if the inequality
\be\label{posdef1} 
\Im \ \tr ( \derl X \derr X\dg  P )  \neq 0 
\ee
holds in $(\xi^0_L,\xi^0_R)$  or if the matrices
\be\label{posdef2} 
\derl X(\xi^0_L,\xi^0_R), \derr X(\xi^0_L,\xi^0_R), X(\xi^0_L,\xi^0_R)
\ee
are linearly independent. Therefore any of the conditions (\ref{posdef1}),(\ref{posdef2}) is a sufficient condition for the existence 
of the surface $\cf$ associated with the solution $X$ of the Euler--Lagrange equations 
(\ref{ELeq}) in the vicinity of the point $(\xi^0_L,\xi^0_R)$. 

Using (\ref{metric}) we can write the formula for Gaussian curvature 
\cite{DoC} 
as
\be
K = \frac{1}{\sqrt{J_{L} J_{R} - G_{LR}^2}} \derr 
\left( \frac{\derl G_{LR} - \frac{1}{2} G_{LR} \derl(\ln J_L)}{\sqrt{J_{L} J_{R} - G_{LR}^2}}\right).
\ee

\section{The Gauss--Weingarten equations}\label{SecGW}

Now we may formally determine a moving frame on the surface $\cf$ and write 
the Gauss--Weingarten equations. 
Let $X$ be a solution of the Euler--Lagrange equations (\ref{ELeq})
 such that ${\rm det}(G)$ is not zero in a neighborhood of a regular point $(\xi_L^0,\xi_R^0)$ in $\Omega$.
Assume also that the surface $\cf$ (\ref{surfaceN}), associated with these equations is described by the moving frame 
$$\vec \tau=(\derl Z, \derr Z, n_{3}, \ldots, n_{N^2-1})^T,$$
where the vectors $\derl Z, \derr Z, n_{3}, \ldots, n_{N^2-1}$ satisfy the normalization conditions
$$ (\derl Z,\derl Z)=J_{L}, \, (\derl Z,\derr Z)= G_{LR}, \, (\derr Z,\derr Z)= J_{R}, $$
\be\label{normnorm} (\derl Z, n_{k})=(\derr Z, n_{k})=0, \, (n_j, n_{k})=\delta_{jk}. \ee
We now show that the moving frame satisfies the Gauss--Weingarten equations
\begin{eqnarray}
\nonumber \derl \derl Z & = & A^L_L \derl Z + A^L_R \derr Z + Q^L_j n_j, \\
\nonumber \derl \derr Z & = & H_j n_j, \\ 
\nonumber \derl n_j & = & \alpha^L_j \derl Z + \beta^L_j \derr Z +s^L_{jk} n_k, \\
\nonumber \derr \derl Z & = & H_j n_j, \\ 
\nonumber \derr \derr Z & = & A^R_L \derl Z + A^R_R \derr z + Q^R_j n_j, \\
\derr n_j & = & \alpha^R_j \derl Z + \beta^R_j \derr Z +s^R_{jk} n_k \label{gweqN},
\end{eqnarray}
where $s^L_{jk}+s^L_{kj}=0$, $s^R_{jk}+s^R_{kj}=0,$ $j,k=3,\ldots, N^2-1$,
$$ \alpha^L_j =  \frac{ H_j G_{LR}- Q^L_j J_{R}}{{\rm det}G}, \, \, \, 
\beta^L_j =  \frac{Q^L_j G_{LR} -   H_j J_{L}}{{\rm det}G},$$ 
$$ \alpha^R_j =  \frac{Q^R_j G_{LR} -   H_j J_{R}}{{\rm det}G}, \, \,  \, 
\beta^R_j =  \frac{ H_j G_{LR}-Q^R_j J_{L}}{{\rm det}G},$$ 
and $A^L_L,A^L_R$ ($A^R_L,A^R_R$ have similar form which can be obtained by exchange $L \leftrightarrow R$) are written as 
\begin{eqnarray}
\nonumber A^L_L  & = & \frac{1}{\det G} \left( J_R (\derl \derl Z,\derl Z) - G_{LR} (\derl \derl Z,\derr Z) \right) \\
A^L_R  & = & \frac{1}{\det G} \left( J_L (\derl \derl Z,\derr Z) - G_{LR} (\derl \derl Z,\derl Z) \right)\label{As}
\end{eqnarray}
where
\begin{eqnarray}
\nonumber (\derl \derl Z,\derl Z)  & = & \frac{1}{2} \tr \left( (\derl \derl X \derl X^\dg + \derl X \derl \derl X\dg )P \right), \\
\nonumber (\derl \derl Z,\derr Z)  & = & -\frac{1}{2} \tr \left( (\derl \derl X \derr X^\dg + \derr X \derl \derl X\dg )P \right. \\
& + & \left. 2 \derl X \derl X\dg (X \derr X\dg + \derr X X\dg) \right).
\end{eqnarray}
Note that in fact we can write it in a compact way
\begin{eqnarray}
\nonumber (\partial_B \partial_B Z,\partial_D Z)  & = & 
(\delta_{B,D}-\frac{1}{2})  \tr \left( (\partial_B \partial_B X \partial_D X^\dg + \partial_D X \partial_B\partial_B X\dg )P \right. \\
& + & \left. 2 \partial_B X \partial_B X\dg (X \partial_D X\dg + \partial_D X X\dg) \right).
\end{eqnarray}

The explicit form of the coefficients $H_j,Q^D_j$ (where $D=L,R$; $j=3,\ldots, N^2-1$) depends on 
the chosen orthonormal basis $\{ n_3, \ldots, n_{N^2-1} \}$ of the normal space to the surface $\cf$ 
at the point $X(\xi_L^0,\xi_R^0)$. Partial information about them will be obtained in 
(\ref{orthddz}).

Indeed, if $\derl Z, \derr Z$ are defined by (\ref{cldx}) for an arbitrary solution $X$ of the Euler--Lagrange equations (\ref{ELeq}),
then by straightforward calculation using (\ref{ELeqmat}) one finds that
\begin{eqnarray}
\nonumber \derl \derr Z  & = &  \derr \derl Z = [\derl P, \derr P] = \\
\nonumber & = & \lambda - \lambda X X\dg - X X\dg \lambda+X (\derl X\dg \derr X - \derr X\dg \derl X) X\dg
 \label{2ndderlr}
\end{eqnarray}
where
$$ \lambda = \derl X \derr X\dg - \derr X \derl X\dg .$$
By computing 
\be\label{11}
\tr \left( \derl \derr Z . \partial_{D} Z \right) = \pm \tr ([\derl P,\derr P].[ \partial_{D} P,P]) = 0, \; D=L,R
\ee
we conclude that $\derl \derr Z$ is perpendicular to the surface $\cf$ and consequently it has the form given in (\ref{gweqN}).

The remaining relations in (\ref{gweqN}) and (\ref{As}) follow as differential consequences 
from the assumed normalizations of the normals (\ref{normnorm}), e.g.
$$ (n_j,n_k)=0, \; j \neq k $$
which gives 
$$0 = (\derl n_j,n_k)+(\derl n_k,n_j) = s^L_{jk} + s^L_{kj}.$$
Similarly
$$ (n_j,\derl Z)=0, \, \, (n_j,\derr Z)=0   $$
by differentiation leads to
$$ (\derr n_j,\derl Z) + (n_j,\derl \derr Z) =0, \, \, (\derr n_j,\derr Z) + (n_j, \derr \derr Z)=0 $$
implying
$$ J_{L} \alpha^R_j + G_{LR} \beta^R_j + H_j =0, \, \, G_{LR} \alpha^R_j + J_{R} \beta^R_j + Q^R_j =0. $$
Consequently, $\alpha^R_j,\beta^R_j$ can be determined in terms of $H_j, Q^R_j$ and of the components of the induced metric 
$G$. The remaining coefficients $\alpha^L_j,\beta^L_j$ are derived in an analogous way by exchanging indices $L \leftrightarrow R$
in the successive differentiations.

The coefficients $A^L_L,\ldots,A^R_R$ are obtained by requiring that $(\partial_{D} \partial_{D} Z - A^D_L \derl Z - A^D_R \derr Z)$
is normal to the surface, i.e. 
\be\label{trdx}
\tr \left( \partial_{B} Z.(\partial_{D} \partial_{D} Z - A^D_L \derl Z - A^D_R \derr Z) \right)=0, \; B,D =L,R.
\ee
From (\ref{mlmr}) and (\ref{cldx}) we find 
\begin{eqnarray}
\nonumber \derl \derl Z & = & [\derl \derl P,P] \\
\nonumber & = & \derl\derl X X\dg -XX\dg \derl\derl X X\dg  +  X\derl \derl X\dg X X\dg \\
& - & X \derl\derl X\dg+2 X \derl X\dg X \derl X\dg-2\derl X X\dg \derl X X\dg. \label{2ndder}
\end{eqnarray}
the expression for $\derr \derr Z$ is obtained by the change of the overall sign and $L \leftrightarrow R $.
After substituting the above expressions into (\ref{trdx}) we solve the resulting linear equations for $A^D_B$.

Let us note that the Gauss--Weingarten equations (\ref{gweqN}) can be written equivalently in the $N \times N$ matrix form
\be
\derl \vec \tau = U \vec \tau, \qquad \derr \vec \tau = V \vec \tau,
\ee
where
\begin{eqnarray}
U & = & \left( \begin{array}{ccccc} A^L_L & A^L_R & Q_3^L & \ldots & Q^L_{N^2-1} \\
0 & 0 & H_3 & \ldots & H_{N^2-1} \\
\alpha^L_3 & \beta^L_3 & s^L_{33} & \ldots & s^L_{3(N^2-1)} \\
\ldots & \ldots & \ldots & \ldots & \ldots \\
\alpha^L_{(N^2-1)} & \beta^L_{(N^2-1)} & s^L_{(N^2-1)3} & \ldots & s^L_{(N^2-1)(N^2-1)} 
\end{array}
\right), \nonumber \\
V & = & \left( \begin{array}{ccccc} 0 & 0 & H_3 & \ldots & H_{N^2-1} \\
A^R_L & A^R_R & Q_3^R & \ldots & Q^R_{N^2-1} \\
\alpha^R_3 & \beta^R_3 & s^R_{33} & \ldots & s^R_{3(N^2-1)} \\
\ldots & \ldots & \ldots & \ldots & \ldots \\
\alpha^R_{(N^2-1)} & \beta^R_{(N^2-1)} & s^R_{(N^2-1)3} & \ldots & s^R_{(N^2-1)(N^2-1)} 
\end{array} \right).
\end{eqnarray}

The Gauss--Codazzi--Ricci equations 
\be\label{gcod1}
\derr U - \derl V + [U,V] =0
\ee
are compatibility conditions for the Gauss--Weingarten equations (\ref{gweqN}).
They are the necessary and sufficient conditions for the local existence of the corresponding surface $\cf$.
It can be easily checked that they are identically satisfied for any solution $X$ of the Euler--Lagrange
equations (\ref{ELeq}).

The second fundamental form and the mean curvature vector of the surface $\cf$ at the regular point $p$ can be expressed, 
according to \cite{Kob,Wil}, as
\bea
\label{IIN1} {\bf II}  =  (\derl \derl Z)^{\perp} {\rm d} \xi_L {\rm d} \xi_L + 2 (\derl \derr Z)^{\perp} {\rm d} \xi_L {\rm d} \xi_R
+(\derr \derr Z)^{\perp} {\rm d} \xi_R {\rm d} \xi_R, & & \\
\label{HN1}
{\bf H}  =  \frac{1}{\det G} \left( J_{R} (\derl \derl Z)^{\perp}  - 2 G_{LR} (\derl \derr Z)^{\perp} 
+ J_{L} (\derr \derr Z)^{\perp} \right), & &
\eea
where $(\;)^{\perp}$ denotes the normal part of the vector.  
In our case the expressions (\ref{IIN1}),(\ref{HN1}) take the form
\begin{eqnarray}
\nonumber {\bf II} & = & \left( Q^L_j {\rm d} \xi_L {\rm d} \xi_L + 2 H_j {\rm d} \xi_L {\rm d} \xi_R
+ Q_j^R {\rm d} \xi_R {\rm d} \xi_R \right) n_j \\
\nonumber & = & (\derl \derl Z-A^L_L \derl Z - A^L_R \derr Z) {\rm d} \xi_L {\rm d} \xi_L 
+ 2 (\derl \derr Z) {\rm d} \xi_L {\rm d} \xi_R +\\
& + &  (\derr \derr Z-A^R_L \derl Z - A^R_R \derr Z) {\rm d} \xi_R {\rm d} \xi_R, \label{IIN2} \\
\nonumber {\bf H} & = & \frac{1}{\det G} \left( J_{R} Q^L_j  - 2 G_{LR} H_j
+ J_{L} Q^R_j \right) n_j \\
\nonumber & = & \frac{1}{\det G} \left( \right. J_{R} (\derl \derl Z-A^L_L \derl Z - A^L_R \derr Z) 
- 2 G_{LR}  (\derl \derr Z)  + \\
& + &  J_{L} (\derr \derr Z-A^R_L \derl Z - A^R_R \derr Z)  \left. \right). \label{HN2}
\end{eqnarray}
The derivatives $\partial_D \partial_B Z$ are expressed explicitly in terms of $X$ in equations (\ref{2ndderlr}) and(\ref{2ndder})
but after substitution of them into (\ref{IIN2}),(\ref{HN2}) get rather complicated, therefore we do not present them here.

\section{The moving frame of a surface in the algebra $\su(N)$}\label{Secmf}

Now we proceed to construct the moving frame of the surface $\cf$ immersed in $\su(N)$ algebra, i.e. matrices 
$ \derl X,\derr X, n_a, \ a=3,\ldots,N^2-1$ satisfying (\ref{normnorm}). 

Let $X$ be a solution of the Euler--Lagrange equations (\ref{ELeq}) and let $(\xi_L^0,\xi_R^0)$ be a regular point 
in $\Omega$, i.e. such that ${\rm det}G(X(\xi_L^0,\xi_R^0)) \neq 0$. Let us denote $X^0=X(\xi_L^0,\xi_R^0)$, 
$Z^0=Z(\xi_L^0,\xi_R^0)$. Taking into account that
$$ \tr(A)= \tr(\Phi A \Phi^{\dagger}), \; \;  A \in \su(N), \  \Phi \in SU(N), $$
we employ the adjoint representation of the group $SU(N)$ in order to bring $ \derl Z,\derr Z, n_a$ to the simplest form possible. 
We shall request
\be\label{Phireq} X^0  =  \Phi \left( \ba{c} \1_{m\times m}  \\
{\bf 0}_{n\times m}   \ea
\right) . \ee
The existence of such $\Phi$ follows from the fact that $G(m,n)$ is a homogeneous space (\ref{Gmn}). In fact $\Phi$ is just any
representant in $SU(N)$ of the equivalence class $[X^0]$ such that
$$ X^0=\left( \ba{ccc} \Phi_{11} & \ldots & \Phi_{1m} \\ 
& \ldots & \\ 
\Phi_{N1} & \ldots & \Phi_{Nm} \ea \right), \; \; \; \;[X^0]=[\Phi] $$
and consequently $\Phi$ is not unique.
Explicit (local) construction of $\Phi$ can be performed algorithmically in such a way that $\Phi$ varies smoothly with smooth change of $X^0$ (for the detailed explanation in the case $G(m,1)=\C P^{m}$ see \cite{Grusno}).

Let us choose an orthonormal basis in $\su(N)$ in the following form
\begin{eqnarray}
\nonumber  (A_{jk})_{ab} & = &  i (\delta_{ja} \delta_{kb} + \delta_{jb} \delta_{ka} ), \; \; 1\leq j<k\leq N, \\
\nonumber  (B_{jk})_{ab} & = &   (\delta_{ja} \delta_{kb} - \delta_{jb} \delta_{ka} ), \; \; 1\leq j<k\leq N, \\
(C_{p})_{ab} & = & i \sqrt{\frac{2}{p(p+1)}} \left( \sum_{d=1}^p \delta_{da} \delta_{db} - p \delta_{p+1,a} \delta_{p+1,b} 
\right), \; \;  1 \leq p \leq N-1.
\end{eqnarray}

Using the definition of $\partial_D Z$  (\ref{cldx}) we find that 
$\Phi^\dagger \partial_D Z(\xi_L^0,\xi_R^0) \Phi$ has the following  block 
structure
\begin{eqnarray}\label{dz0}
 \partial_D^\Phi Z^0 & \equiv & \Phi^\dagger \partial_D Z(\xi_L^0,\xi_R^0) \Phi  =    
\left( \ba{cc} {\bf 0}_{m\times m} & - \partial_{D}^{\Phi} Z^\dg \\ 
\partial_{D}^{\Phi} Z & {\bf 0}_{n\times n} \ea \right)
\end{eqnarray}
(where $\partial_{D}^{\Phi} Z$ are defined by (\ref{dz0})).

When $\Phi$ satisfying (\ref{Phireq}) is found, the construction of the moving frame can proceed as follows.
Assume that one finds, using a variant of Gramm-Schmidt orthogonalization procedure, the orthonormal vectors 
$$ \tilde A_{aj}, \tilde B_{aj}, \; a = 1,\ldots,m, \; j=m+1,\ldots, N, \; a+j>m+2 $$
satisfying
$$ (\partial_{D}^\Phi Z^0, \tilde A_{aj}) =0, \, (\partial_{D}^\Phi Z^0, \tilde B_{bj}) =0 $$
and 
\bea
\nonumber {\rm span} (\partial_{D}^\Phi Z^0,\tilde A_{aj}, \tilde B_{aj})_{D=L,R, a = 1,\ldots,m, \; j=m+2,\ldots, N, \; a+j>m+2}  & = & \\ 
{\rm span} (A_{aj},B_{aj})_{a = 1,\ldots,m, \; j=m+1,\ldots, N}. & &
\label{GrammSch}
\eea
We identify the remaining tilded and untilded matrices
$$ \tilde A_{jk} = A_{jk}, \ \tilde B_{jk} = B_{jk}, \ \tilde C_{p} = C_{p}, $$
where $a,j = 1,\ldots,m$ or $a,j = m+1,\ldots,N$, $1 \leq p \leq N-1.$
As a result, from Gramm-Schmidt orthogonalization and (\ref{dz0}) we get
$$ (\partial_{D}^\Phi Z^0,\tilde A_{ak})=(\partial_{D}^\Phi Z^0,\tilde B_{ak})=(\partial_{D}^\Phi Z^0,\tilde C_{p})=0 $$
and 
$$ (\tilde A_{ai},\tilde A_{bk}) = (\tilde B_{ai},\tilde A_{bk}) = \delta_{ab} \delta_{ik}, (\tilde C_{p},\tilde C_{q})=\delta_{pq},$$
$$ (\tilde A_{1i},\tilde B_{jk}) = (\tilde A_{1i},\tilde C_{p}) = (\tilde B_{1i},\tilde C_{p}) = 0  $$
(where indices run through all the values for which the respective tilded matrices are defined).

Therefore, under the above given assumptions and notation, we can state the following 
\begin{prop}
Let the moving frame of the surface $\cf$ in the neighborhood $\Upsilon$ of point $Z^0=Z(\xi_L^0,\xi_R^0)$ be
\begin{eqnarray}
\nonumber \derl Z  & = & \Phi \derl^\Phi Z \Phi^\dagger, \\ 
\nonumber \derr Z  & = & \Phi \derr^\Phi Z \Phi^\dagger, \\
\nonumber n^A_{jk} & = & \Phi \tilde A_{jk} \Phi^\dagger, \\
\nonumber n^B_{jk} & = & \Phi \tilde B_{jk} \Phi^\dagger,  \\
 n^C_{p} & = & \Phi \tilde C_{p} \Phi^\dagger. \label{movframeN}
\end{eqnarray}
where indices run through the values for which $\tilde A,\tilde B,\tilde C$ are defined and $ \Phi(\xi_L,\xi_R)$ on $\Upsilon$ satisfies
$$ X(\xi_L,\xi_R) = \Phi(\xi_L,\xi_R) \left( \ba{c} \1_{m\times m}  \\
{\bf 0}_{n\times m}   \ea
\right) .$$
Then (\ref{movframeN}) satisfies the normalization conditions (\ref{normnorm})
and consequently the Gauss--Weingarten equations (\ref{gweqN}).
\end{prop}
\smallskip

Note that the first two lines of (\ref{movframeN}) are equivalent to (\ref{dz0}). 
The remaining lines of (\ref{movframeN}) give a rather explicit 
description of normals to the surface $\cf$. Since the construction is local, 
we don't have an apriori control of the orientation of normals. Of course, in the 
neighborhood where the procedure is applied the normals have the same orientation.

The explicit form of the moving frame (\ref{movframeN}) might be quite complicated
 because of the orthogonalization process involved in the construction of 
$$ n^A_{aj}, n^B_{aj}, \; a = 1,\ldots,m, \; j=m+1,\ldots, N, \; a+j>m+2  $$
(i.e. in the construction of $\tilde A_{1j}, \tilde B_{1j}$).
On the other hand, the remaining normals 
$$  n^A_{ak}, \  n^B_{ak}, \  n^C_{p} $$
where $ a,j = 1,\ldots,m$ or $a,j = m+1,\ldots,N$, $1 \leq p \leq N-1$
can be constructed immediately after finding $\Phi$.\footnote{
In fact, one particular combination of normals $n^C_{p}$ of the form
$$ n_P = i \sqrt{2} \left( \sqrt{\frac{N-m}{mN}} \ \1 - \sqrt{\frac{N-m}{m(N-m)}} \ P \right) $$
can be constructed from $X$, i.e. $P$, alone, without the knowledge of corresponding $\Phi$.
}

If we choose other group element $\Phi$ satisfying (\ref{Phireq}),
the constructed normals would have been rotated by a local (gauge) transformation from the subgroup of $SU(N)$ leaving 
$\derl Z,\derr Z$ invariant.

It is worth noting that from the equations (\ref{2ndderlr}),(\ref{2ndder}) immediately follows that 
\begin{eqnarray}
\label{orthddz} (\derl \derl Z)^\perp, (\derr \derr Z)^\perp  & \in & \  {\rm span} (n^A_{1j},n^B_{1j})_{a = 1,\ldots,m, j=m+1,\ldots, N}, \\
\nonumber (\derl \derr Z)^\perp=\derl \derr Z & 
\in &  \  {\rm span} (n^A_{jk},n^B_{jk},n^C_{p})_{a,j = 1,\ldots,m, \ {\rm or} \ a,j=m+1,\ldots, N,p<N},
\end{eqnarray}
i.e. $\derl \derr Z$ is orthogonal to $\derl \derl Z,\derr \derr Z$ (and also to the surface, see (\ref{gweqN})).

\section{Final remarks}\label{concl}

The technique presented above for finding surfaces associated with complex Grassmannian sigma models defined 
on Minkowski space can be seen to give a  rather detailed analytical description of the surfaces in question.
This description provides  effective tools for constructing surfaces without 
reference to additional considerations, proceeding directly from the given complex Grassmannian sigma model equations (\ref{ELeq}). 
Through the use of Cartan's language of moving frames 
we derived via this sigma model the structural equations of two--dimensional smooth surfaces immersed in $\su(N)$ algebra.
It allows to find the first 
and second fundamental forms of the surfaces as well as the relations between them as expressed in the Gauss--Weingarten and 
Gauss--Codazzi--Ricci equations. 
An extension of the classical Enneper--Weierstrass representation of surfaces in multidimensional spaces,
expressed in terms of any nonsingular solution of (\ref{ELeq}), was presenetd in explicit form.

\subsection*{Acknowledgments}

This work was supported in part by research grants from NSERC of Canada.
Libor \v Snobl acknowledges a postdoctoral fellowship awarded by the Laboratory of Mathematical Physics of the CRM, 
Universit\'e de Montr\'eal. The authors thank Pavel Winternitz for helpful and interesting discussions on the topic
of this paper.


\begin{thebibliography}{99}
\bibliographystyle{plain}

\bibitem{Abl}
Ablowitz, M., Chakravarty, S. and  Halburd, R., 
{\em J. Math. Phys.}  \textbf{44}, 3147-3173 (2003).

\bibitem{David}
David, F.,  Ginsparg, P. and  Zinn-Justin, Y. (editors), 
{\em Fluctuating Geometries in Statistical Mechanics and Field Theory}, Elsevier (Amsterdam,
1996).

\bibitem{Nel}
Nelson, D.,  Piran, T. and  Weinberg, S., {\em Statistical Mechanics of Membranes and Surfaces}, World Scientific (Singapore, 1992).

\bibitem{Cha}
Charvolin, J., Joanny, J.F.,  Zinn-Justin, J., {\em Liquids at Interfaces}, Elsevier (Amsterdam, 1989).

\bibitem{Ou}
Ou-Yang, Z.,  Lui, J. and  Xie, Y., {\em Geometric Methods in Elastic Theory of Membranes in Liquid Crystal Phases}, 
World Scientific (Singapore, 1999).

\bibitem{Saf}
Safram, S.A., {\em Statistical Thermodynamics of Surfaces, Interfaces and Membranes}, Addison-Wesley (New York, 1994).

\bibitem{Grusno}
A. M. Grundland and L. \v Snobl, Description of surfaces associated with $\C P^(N-1)$ sigma models on Minkowski space, 
submitted for publication, [math.DG/0405513].

\bibitem{Grusnosym}
Grundland, A.M. and \v Snobl, L., Surfaces in $\su(N)$ algebra via $\C P^{N-1}$ sigma models on Minkowski space, 
{\em Proceedings of  XI International Conference on Symmetry Methods in Physics}, (Prague 2004), in press.

\bibitem{Sas}
Sasaki, R.,
{\em Phys. Lett} \textbf{130B}, 69--72 (1983).

\bibitem{Din}
Din, A.M., Nonlinear technique in two dimensional Grassmannian sigma models, Lecture Notes Math. 1139 (1983) 253-262.

\bibitem{Gro}
Gross, D.J., Piran, T. and  Weinberg, S.,
{\em Two--dimensional quantum gravity and random surfaces}, World Scientific (Singapore, 1992). 

\bibitem{Bar}
Barrett, J., Gibbons, G.W., Perry, M.J., and Ruback, P., 
{\em Int. J. Mod. Phys.} \textbf{A 9} 
1457-1493 (1994).

\bibitem{Bre}
Brezin, E., Izykson, C., Zinn--Justin, J. and Zuber, J.B.,
{\em Phys. Lett.} \textbf{82B} 442 (1979).

\bibitem{Lus}
Luscher, M.  and Pohlmeyer, K., {\em Nucl.Phys.} \textbf{B137} 46 (1978).

\bibitem{Koi}
Koikawa, T. and Sasaki, R., {\em Phys. Lett.} \textbf{B124} 85 (1983).

\bibitem{Pie}
Piette, B., 
{\em J. Math. Phys.}  \textbf{29}, 2190--2196 (1988).

\bibitem{Jac}
Jacques, M. and Saint-Aubin, Y.,
{\em  J. Math. Phys.} \textbf{28}, 2463--2479 (1987). 

\bibitem{DoC}
do Carmo, M. P., {\em Riemannian Geometry}, Birkh\"auser
(Boston, 1992).

\bibitem{Kob}
Kobayashi S. and  Nomizu, K., {\em Foundation of Differential Geometry}, John Wiley (New York, 1963).

\bibitem{Wil}
Willmore, T.J., {\em Riemannian Geometry}, Clarendon (Oxford, 1993).


\end{thebibliography}
\end{document}